\documentclass[11pt]{article}
\usepackage{amsmath,amsfonts,latexsym,amssymb,chicago,graphicx,annals,here}
\advance \topmargin by -\headheight \advance \topmargin by
-\headsep

\newcommand{\mycite}[1]{{\small \sc \citeNP{#1}}}

\def\R{\mathbb R}

\def\ai{\mbox{Ai}}
\def\bi{\mbox{Bi}}

\def\f{\phi}
\def\l{\lambda}
\def\eop{\hfill\mbox{$\Box$}\newline}

\newtheorem{corollary}{Corollary}[section]
\newtheorem{lemma}{Lemma}[section]

\begin{document}
\title{The maximum of Brownian motion minus a parabola}
\author{Piet Groeneboom}
\date{Submitted: 10-5-2010. Accepted: 10-5-2010}
\affiliation{Delft University of Technology, Mekelweg 4, 2628 CD Delft, The Netherlands,
p.groeneboom@tudelft.nl; http://dutiosc.twi.tudelft.nl/~pietg/}
\AMSsubject{60J65,60J75}
\maketitle
\begin{abstract}
We derive a simple integral representation for the distribution of the maximum of Brownian motion minus a parabola, which can be used for computing the density and moments of the distribution, both for one-sided and two-sided Brownian motion.
\end{abstract}

\section{Introduction}
\label{intro}
It is the purpose of this note to show how one can easily obtain information on properties of the distribution of the maximum of Brownian motion minus a parabola from \mycite{piet:89}. In fact, Corollary 3.1 in that paper gives the joint distribution of both the maximum and the location of the maximum. In the latter paper most attention is on the distribution of the {\it location} of the maximum, which is derived from this corollary. The reason for the emphasis on the distribution of the location of the maximum is that this distribution very often occurs as limit distribution in the context of isotonic regression; one could say that it is a kind of ``normal distribution" in that context. But one can of course also derive the distribution of the maximum itself from this corollary and at the same time deduce  numerical information, as will be shown below.

Numerical information on the density, quantiles and moments of the location of the maximum is given in \mycite{piet_jon:01}, which in turn relies on section 4 of \mycite{piet:85}.

\section{Representations of the distribution of the maximum}
\label{characterization}
Let $F_c$ be the distribution function of the maximum of $W(t)-ct^2,\, t\ge0,$ where $W$ is one-sided Brownian motion (in standard
scale and without drift). Then, according to Theorem 3.1 of \mycite{piet:89}, $F_c$ has the representation
\begin{equation}
\label{F_repr}
F_c(x)={\psi}_{x,c}(0),
\end{equation}
where the function $\psi_{x,c}:\R\to\R_+$ has Fourier transform
\begin{equation}
\label{psi_hat}
\hat\psi_{x,c}(\l)=\int_{-\infty}^{\infty}e^{i\l
s}\psi_{x,c}(s)\,ds=\frac{\pi\left\{\ai(i\xi)\bi(i\xi+z)-\bi(i\xi)\ai(i\xi+z)\right\}}{(2c^2)^{1/3}\ai(i\xi)}\,\,,
\end{equation}
and where $\xi=\left(2c^2\right)^{-1/3}\l,$ and $z=(4c)^{1/3}x$. It follows that the corresponding density $f_c$ has the
representation
\begin{equation}
\label{f_repr}
f_c(x)=\f_{x,c}(0)\stackrel{\mbox{\mbox{\small def}}}=\frac{\partial}{\partial x}{\psi}_{x,c}(0),
\end{equation}
where $\f_{x,c}$ has Fourier transform
\begin{equation}
\label{phi_hat}
\hat\f_{x,c}(\l)=\int_{-\infty}^{\infty}e^{i\l s}\f_{x,c}(s)\,ds
=\frac{\pi(4c)^{1/3}\left\{\ai(i\xi)\bi'(i\xi+z)-\bi(i\xi)\ai'(i\xi+z)\right\}}{(2c^2)^{1/3}\ai(i\xi)}\,\,,
\end{equation}

The (symmetric) density $g_c$ of the maximum $M$ of $W(t)-ct^2,\,t\in\R$, where $W$ is {\it two-sided} Brownian motion, originating from zero, therefore
has the representation
\begin{equation}
\label{phi_hat2}
g_c(x)=2f_c(x)F_c(x)=2\f_{x,c}(0)\psi_{x,c}(0)=\frac1{2\pi^2}\int_{-\infty}^{\infty}\hat\psi_{x,c}(u)\,du\int_{-\infty}^{\infty}\hat\f_{x,c}(u)\,du,\,x>0,
\end{equation}
since the distribution function of $M$ is the maximum of the two maxima one gets to the right and to the left of zero.
Note that these two maxima are independent, since two-sided Brownian motion is started independently to the right and to the
left, starting at zero. It is also obvious that these two maxima have the same distribution. Interestingly, the situation is
more complicated for the {\it location} of the maximum!

Note that this gives the complete characterization of the distribution of the maximum of Brownian motion with parabolic drift. The purpose of this note, however, is to show how one can deduce useful numerical information from this.

The two fundamental solutions of the Airy differential equation are $\ai$ and $\bi$ which are unbounded on different regions of the complex plane. For the purpose of computing moments, etc., it is easier to only work with the solution $\ai$, so we want to get rid of $\bi$. To this end we simply use Cauchy's formula.

We have the following lemma.

\begin{lemma}
\label{df_N}
Let $N_c$ be defined by
$$
N_c=\max_{x\ge0}\left\{W(x)-cx^2\right\},
$$
So $N_c$ is the maximum for the one-sided case. Then the distribution function $F_c=F_{N_c}$ of $N_c$ is given by:
$$
F_c(x)=1-\int_{(4c)^{1/3}x}^{\infty}{\mathrm {Ai}}(u)\,du -2\,{\mathrm {Re}}\left\{e^{-i\pi/6}\int_0^{\infty}\frac{{\mathrm {Ai}}\left(e^{-i\pi/6}u\right){\mathrm {Ai}}(iu+(4c)^{1/3}x)}{{\mathrm {Ai}}(iu)}\,du\right\},\,x>0.
$$
\end{lemma}

One can use this representation to compute the distribution function $F_c$ in one line in, for example, Mathematica, and the result of this computation is shown below in Figure \ref{fig:df}, where we take $c=1/2$.

\begin{figure}[!ht]
\begin{center}
\includegraphics{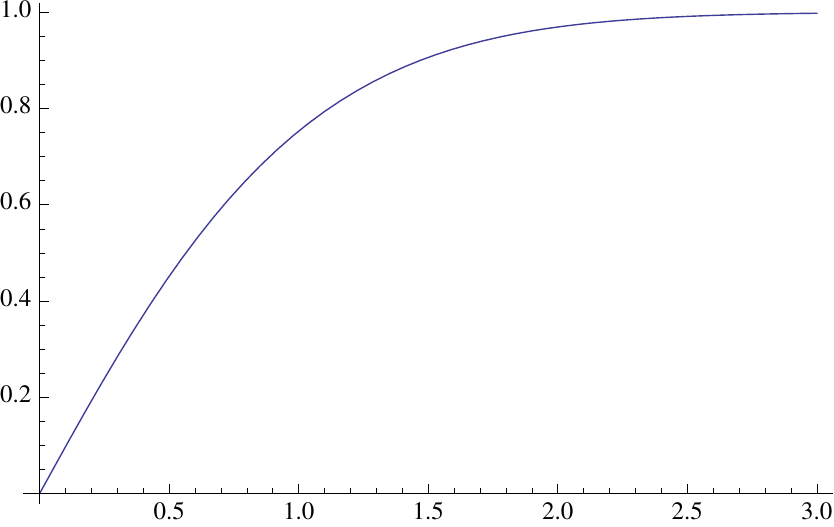}
\end{center}
\caption{The distribution function $F_{N_c}$, for $c=1/2$.}
\label{fig:df}
\end{figure}

\vspace{0.3cm}
\noindent
{\bf Proof of Lemma \ref{df_N}.}
After the change of variables $u=(2c^2)^{-1/3}\xi$ we
get for the corresponding distribution function $F_c$, still taking $z=(4c)^{1/3}x$,
\begin{equation}
\label{F_repr2}
F_c(x)=\frac12\int_{-\infty}^{\infty}\frac{\left\{\ai(iu)\bi(iu+z)
-\bi(iu)\ai(iu+z)\right\}}{\ai(iu)}\,du.
\end{equation}
Using
$$
\bi(z)=i\ai(z)-2i e^{\pi i/3}\ai\bigl(z e^{-2\pi i/3}\bigr),
$$
which is 10.4.9 in \mycite{AbSte:64}, we can write:
\begin{eqnarray}
\label{F_repr3}
&&\frac12\int_0^{\infty}\frac{\left\{\ai(iu)\bi(iu+z)
-\bi(iu)\ai(iu+z)\right\}}{\ai(iu)}\,du\nonumber\\
&&=e^{-i\pi/6}\int_0^{\infty}\ai\left(e^{-i\pi/6}u+ze^{-2i\pi/3}\right)\,du
-e^{-i\pi/6}\int_0^{\infty}\frac{\ai\left(e^{-i\pi/6}u\right)\ai(iu+z)}
{\ai(iu)}\,du.
\end{eqnarray}
By Cauchy's formula we can reduce the first integral on the right-hand side of (\ref{F_repr3}) to:
\begin{equation}
\label{integral1}
\int_0^{\infty}\ai\left(u+ze^{-2i\pi/3}\right)\,du.
\end{equation}
Differentiation w.r.t.\ $z$ yields:
$$
e^{-2i\pi/3}\int_0^{\infty}\ai'\left(u+ze^{-2i\pi/3}\right)\,du=-e^{-2i\pi/3}\ai\left(ze^{-2i\pi/3}\right).
$$

We similarly have:
\begin{eqnarray}
\label{F_repr4}
&&\frac12\int_{-\infty}^0\frac{\left\{\ai(iu)\bi(iu+z)
-\bi(iu)\ai(iu+z)\right\}}{\ai(iu)}\,du\nonumber\\
&&=e^{i\pi/6}\int_0^{\infty}\ai\left(e^{i\pi/6}u+ze^{2i\pi/3}\right)\,du
-e^{i\pi/6}\int_0^{\infty}\frac{\ai\left(e^{i\pi/6}u\right)\ai(-iu+z)}
{\ai(-iu)}\,du.
\end{eqnarray}
Again using Cauchy's formula we can reduce the first integral on the right-hand side of (\ref{F_repr4}) to:
\begin{equation}
\label{integral2}
\int_0^{\infty}\ai\left(u+ze^{2i\pi/3}\right)\,du.
\end{equation}
Differentiation w.r.t.\ $z$ yields:
$$
e^{2i\pi/3}\int_0^{\infty}\ai'\left(u+ze^{2i\pi/3}\right)\,du=-e^{2i\pi/3}\ai\left(ze^{2i\pi/3}\right).
$$
So the derivative w.r.t.\ $z$ of the sum of the two integrals (\ref{integral1}) and (\ref{integral2}) is given by
\begin{equation}
\label{diff}
-e^{-2i\pi/3}\ai\left(ze^{-2i\pi/3}\right)-e^{2i\pi/3}\ai\left(ze^{2i\pi/3}\right)
=\ai(z).
\end{equation}
For the latter relation, see (10.4.7) in \mycite{AbSte:64}.

So we get:
\begin{equation}
\label{integral3}
\int_0^{\infty}\ai\left(u+ze^{-2i\pi/3}\right)\,du
+\int_0^{\infty}\ai\left(u+ze^{2i\pi/3}\right)\,du=\int_0^z\ai(u)\,du+k,
\end{equation}
for some constant $k$. For $z=0$ we get:
$$
\int_0^{\infty}\ai(u)\,du+\int_0^{\infty}\ai(u)\,du=2\int_0^{\infty}\ai(u)\,du=k=2/3.
$$
Hence $k=2/3$ and
\begin{eqnarray}
\label{integral4}
&&\int_0^{\infty}\ai\left(u+ze^{-2i\pi/3}\right)\,du
+\int_0^{\infty}\ai\left(u+ze^{2i\pi/3}\right)\,du=\int_0^z\ai(u)\,du+2/3\nonumber\\
&&=1-\int_z^{\infty}\ai(u)\,du.
\end{eqnarray}
Note that this implies:
\begin{equation}
\label{lim_integral}
\lim_{z\to\infty}\biggl\{\int_0^{\infty}\ai\left(u+ze^{-2i\pi/3}\right)\,du
+\int_0^{\infty}\ai\left(u+ze^{2i\pi/3}\right)\,du\biggr\}=1.
\end{equation}

The second term is given by
\begin{eqnarray}
\label{sec_term}
&&-e^{-i\pi/6}\int_0^{\infty}\frac{\ai\left(e^{-i\pi/6}u\right)\ai(iu+z)}
{\ai(iu)}\,du
-e^{i\pi/6}\int_0^{\infty}\frac{\ai\left(e^{i\pi/6}u\right)\ai(-iu+z)}
{\ai(-iu)}\,du\nonumber\\
&&=-2\mbox{Re}\left\{e^{-i\pi/6}\int_0^{\infty}\frac{\ai\left(e^{-i\pi/6}u\right)\ai(iu+z)}
{\ai(iu)}\,du\right\}.
\end{eqnarray}
\eop

\begin{corollary}
\label{dens_N}
The density of $N_c$ is given by:
$$
f_c(x)=(4c)^{1/3}\left\{{\mathrm {Ai}}\bigl((4c)^{1/3}x\bigr) -2\,\mbox{\rm Re}\left(e^{-i\pi/6}\int_0^{\infty}\frac{{\rm {\mathrm {Ai}}}\left(e^{-i\pi/6}u\right){\mathrm {Ai}}\,'\bigl(iu+(4c)^{1/3}x\bigr)}{\mathrm {Ai}(iu)}\,du\right)\right\},\,x>0.
$$
\end{corollary}

\noindent
{\bf Proof.} This follows by straightforward differentiation from Lemma \ref{df_N}.\eop

\begin{figure}[!ht]
\begin{center}
\includegraphics{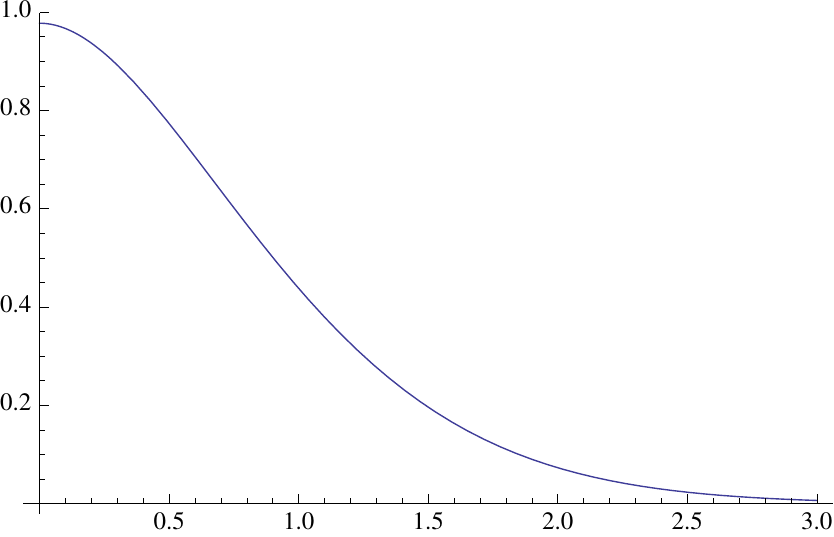}
\end{center}
\caption{The density $f_{N_c}$, for $c=1/2$.}
\label{fig:density}
\end{figure}

\vspace{0.3cm}
\noindent
A picture of the density $f_{N_c}$, for $c=1/2$, is given in Figure \ref{fig:density}.
By the remarks above, we also have:

\begin{corollary}
\label{dens_M}
The density $g_c$ of the maximum $M$ of $W(t)-ct^2,\,t\in\R$, where $W$ is {\it two-sided} Brownian motion, originating from zero, if given by
$$
g_c(x)=2f_c(x)F_c(x),\,x>0,
$$
where $f_c$ is given by Corollary \ref{dens_N} and $F_c$ by Lemma \ref{df_N}.
\end{corollary}

\begin{figure}[!ht]
\begin{center}
\includegraphics{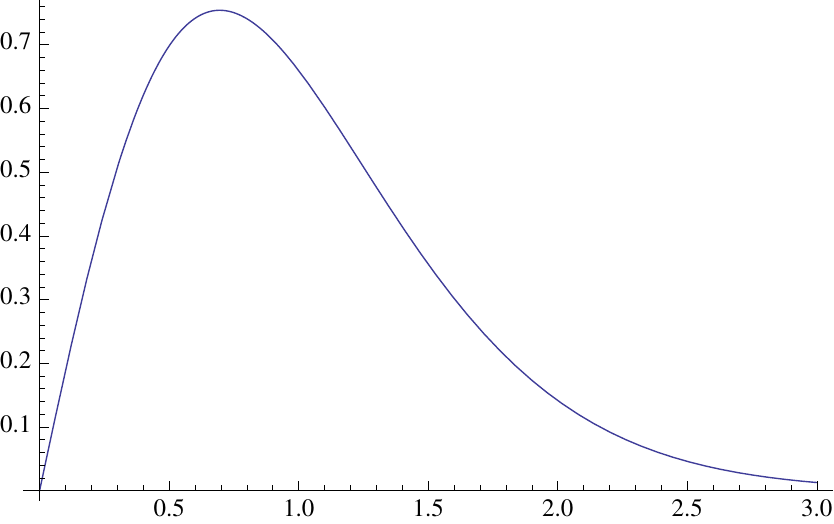}
\end{center}
\caption{The density $f_{M_c}=g_c$ of the maximum for two-sided Brownian motion and $c=1/2$.}
\label{fig:2sided_density}
\end{figure}

\vspace{0.3cm}
\noindent
A picture of the density $f_{M_c}$, for $c=1/2$ and two-sided Brownian motion, is given in Figure \ref{fig:2sided_density}.

\section{Concluding remarks}
\label{section:conclusion}
The densities of the maximum and location of the maximum of Brownian motion minus a parabola were originally studied by solving partial differential equations. For example, if we denote the location of the maximum of two-sided Brownian motion minus the parabola $y=t^2$ by $Z$, then the density of $Z$ is expressed in \mycite{chernoff:64} in terms of the solution of the heat equation
$$
\frac{\partial}{\partial t}u(t,x)=-\tfrac12\frac{\partial^2}{\partial x^2}u(t,x),
$$
for $x\le t^2$, under the boundary conditions
$$
u(t,t^2)\stackrel{\mbox{def}}=\lim_{x\uparrow t^2}u(t,x)=1,\qquad \lim_{x\downarrow -\infty}u(t,x)=0,\qquad t\in\R.
$$
If $u(t,x)$ is the (smooth) solution of this equation, the density $f_Z$ of $Z$ is given by
$$
f_Z(t)=\tfrac12 u_2(-t)u_2(t),\,x\in\R,
$$
where (as in \mycite{piet:85}) the function $u_2$ is defined by
$$
u_2(t)=\lim_{x\uparrow t^2}\frac{\partial}{\partial x}u(t,x).
$$
The original computations of this density were indeed based on numerically solving this partial differential equation (as I learned from personal communications by Herman Chernoff and Willem van Zwet). However, it is very hard to solve this equation numerically sufficiently accurately for negative values of $t$, since we have, by (4.25) in \mycite{piet:85}:
$$
u_2(t)\sim c_1\exp\left\{-\tfrac23|t|^3-c|t|\right\},\,t\to-\infty,
$$
where $c\approx 2.9458\dots$ and $c_1\approx 2.2638\dots$.

\begin{figure}[!ht]
\begin{center}
\includegraphics{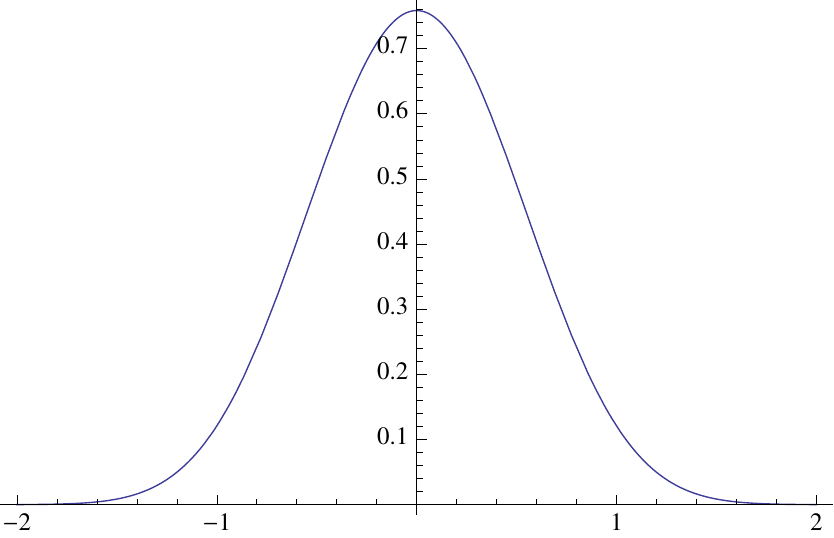}
\end{center}
\caption{The density $f_{Z}$ of the location of the maximum of $W(t)-t^2,\,t\in\R$.}
\label{fig:loc_max}
\end{figure}

At present the situation is drastically different, since we have much more analytical information about the solution and, moreover, can use advanced computer algebra packages. One only needs one line in Mathematica to compute the density $f_Z$, since, by (3.8) in \mycite{piet:89}, $f_Z$ is given by $f_Z(x)=\tfrac12 \f(x)\f(-x)=\tfrac12 u_2(x)u_2(-x)$, where:
$$
\f(x)=\frac1{2^{2/3}\pi}\int_{-\infty}^{\infty}\frac{e^{-iux}}{\ai(i2^{-1/3}u)}\,du,
$$
and where one can even allow the boundaries $-\infty$ and $\infty$ in the numerical integration (in Mathematica).
A picture of the density $f_Z$, obtained from just using this definition in Mathematica, is given in Figure 
\ref{fig:loc_max}. 

However, if one wants to get very precise information about the tail behavior of the density or the behavior close to zero, it is better to use power series expansions or asymptotic expansions, which are different in a neighborhood of zero from the representation for large values of the argument. Details on this are given in \mycite{piet:85} and \mycite{piet_jon:01}.

More details on the history of the subject are given in \mycite{perman_wellner96} and \mycite{janson:10}.
In the latter manuscript also more details on the distribution of the maximum of Brownian motion minus a parabola are given. Their results seem to be in complete agreement with some numerical computations, based on the representations given in section \ref{characterization}.

\vspace{0.3cm}
\noindent
{\bf Acknowledgement}
I want to thank Neil O'Connell for inviting me to submit this note to the Electronic Journal of Probability.

\end{document}